\documentclass[12pt]{article}

\usepackage{graphicx}
\usepackage{enumerate}
\usepackage{amssymb}
\usepackage{amsmath}
\usepackage{mathrsfs}
\usepackage{latexsym}
\usepackage{psfrag}
\usepackage{mathrsfs}
\usepackage{verbatim}
\usepackage{xcolor}
\usepackage{comment}

\usepackage{enumerate,lineno,setspace,float}

\usepackage[utf8]{inputenc}
\usepackage{amsfonts, amsmath, amssymb, amsthm}
\usepackage{graphicx}
\usepackage{comment}
\usepackage{tikz, pgffor}
\usepackage{hyperref}
\usepackage[numbers,sort&compress]{natbib}
\usepackage[capitalize]{cleveref}
\hypersetup{
    colorlinks=true,       
    linkcolor=blue,        
    citecolor=blue,        
}

\allowdisplaybreaks
\newtheorem{theorem}{Theorem}

\newtheorem{corollary}[theorem]{Corollary}
\newtheorem{definition}{Definition}

\newtheorem{conjecture}{Conjecture}

\newcommand{\Href}[1]{\hyperref[#1]{\Cref{#1}}}
\renewcommand{\href}[1]{\hyperref[#1]{\ref{#1}}}

\def\leq{\leqslant}
\def\geq{\geqslant}

\def\deg{{\rm deg}}

\usepackage{color}
\definecolor{VeryLightBlue}{rgb}{0.9,0.9,1}
\definecolor{LightBlue}{rgb}{0.8,0.8,1}
\definecolor{MidBlue}{rgb}{0.3,0.3,1}
\definecolor{DarkBlue}{rgb}{0,0,0.6}
\definecolor{Blue}{rgb}{0,0,1}
\definecolor{Gold}{rgb}{1,0.843,0}
\definecolor{LightGreen}{rgb}{0.88,1,0.88}
\definecolor{MidGreen}{rgb}{0.6,1,0.6}
\definecolor{DarkGreen}{rgb}{0,0.6,0}
\definecolor{VeryLightYellow}{rgb}{1,1,0.9}
\definecolor{LightYellow}{rgb}{1,1,0.6}
\definecolor{MidYellow}{rgb}{1,1,0.5}
\definecolor{DarkYellow}{rgb}{1,1,0.01}
\definecolor{DarkPurple}{rgb}{.6,0,1}
\definecolor{Red}{rgb}{1,0,0}
\definecolor{VeryLightRed}{rgb}{1,0.9,0.9}
\definecolor{LightRed}{rgb}{1,0.8,0.8}
\definecolor{MidRed}{rgb}{1,0.53,0.55}

\newcommand{\darkpurple}[1]{{\color{DarkPurple}{#1}}}
\newcommand{\li}[1]{\darkpurple{\bf [Li: #1]}}

\title{On $k$-Shifted Antimagic Spider Forests} 
\author{Fei-Huang Chang
\thanks{Division of Preparatory Programs for Overseas Chinese Students, National Taiwan Normal
University, New Taipei City 24449, Taiwan. 
{\tt Email:cfh@ntnu.edu.tw.} Supported by MOST-109-2115-M-003-005.}
\and
Wei-Tian Li
\thanks{Department of Applied Mathematics, National Chung Hsing University, Taichung 40227, Taiwan. 
{\tt Email:weitianli@nchu.edu.tw.} Supported by MOST-110-2115-M-005 -005 -MY2.}
\and
Daphne Der-Fen Liu
\thanks{California State University, Los Angeles, USA.  
{\tt Email:dliu@exchange.calstatela.edu}.  Partially supported by NSF grant DM 1600778 and Cal State LA Provost Research Fellow grant.}
\and
Zhishi Pan
\thanks{Department of Mathematics, Tamkang University, New Taipei City,  251301, 
Taiwan. 
{\tt Email:zhishi@mail.tku.edu.tw.} Supported by MOST-109-2115-M-032-003.
}
}

\begin{document}

\maketitle

\baselineskip=18pt

   



\begin{abstract}
Let $G(V,E)$ be a simple graph with $m$ edges. For a given integer $k$, a $k$-shifted antimagic labeling is a bijection $f: E(G) \to \{k+1, k+2, \ldots, k+m\}$ such that all vertices have different vertex-sums, where the vertex-sum of a vertex $v$ is the total of the labels assigned to the edges incident to $v$. 
A graph $G$ is {\it $k$-shifted antimagic} if it admits a $k$-shifted antimagic labeling. For the special case when $k=0$, a  $0$-shifted antimagic labeling is known as {\it antimagic labeling}; and $G$ is {\it antimagic} if it admits an antimagic labeling. A spider is a tree with exactly one vertex of degree greater than two. A spider forest is a graph where each component is a spider. In this article, we prove that certain spider forests are  $k$-shifted antimagic for all $k \geq 0$. In addition, we show that for a  spider forest $G$ with $m$ edges, there exists a positive integer $k_0< m$ such that $G$ is $k$-shifted antimagic for all $k \geq k_0$ and $k \leq -(m+k_0+1)$.
\end{abstract}

{\bf keywords:} Antimagic labeling; $k$-shifted antimagic labeling; spider forest. 
\section{Introduction}


All graphs considered in this paper are finite, simple, and undirected. For any two integers $a \leq  b$, denote $[a,b] = \{x\in\mathbb{Z}: a \leq  x \leq  b\}$. In 1990, Hartsfield and Ringel~\cite{Ringel} introduced the following labeling problem on graphs:  

\begin{definition}\rm
\label{antimagic}
Let $G(V,E)$ be a graph with $m$ edges. For a bijection $f: E  \to [1, m]$,  define the {\it vertex-sum} of $v \in V(G)$, $\phi_f(v)$, as the total of the labels of edges incident to $v$. That is, $\phi_f(v) = \sum_{e \cap v \neq \emptyset} f(e)$.  The {\it antimagic condition} is defined as
$$
\phi_f(v) \neq \phi_f(u), \ \ \ \mbox{if} \ u \neq v. 
$$
If the antimagic condition is satisfied, then $f$ is called an {\it antimagic labeling} of $G$. Further,  $G$ is {\it antimagic} if $G$ admits an antimagic labeling. When the bijective mapping $f$ is clear in the context, we simply denote $\phi_f(v)$ by $\phi(v)$. 
\end{definition}

Clearly, the 2-path $P_2$ (one-edge path) is not antimagic. For connected graphs with at least three vertices, Hartsfield and Ringel \cite{Ringel} proved that some families of graphs, such as paths, stars, cycles, and complete graphs are antimagic, and proposed the following two conjectures:  

\begin{conjecture} 
{\rm (\cite{Ringel})} 
\label{connected antimagic}
Every connected graph on at least three vertices is antimagic.
\end{conjecture}

\begin{conjecture}
{\rm (\cite{Ringel})} 
\label{tree antimagic}
Every tree on at least three vertices is antimagic.
\end{conjecture}

\Href{connected antimagic} and \Href{tree antimagic} have been studied intensively.  
Many classes of graphs have been proved to be antimagic~\cite{dense,regular,regular bipartite}. For a comprehensive list of known  results, we refer the readers to the survey by Gallian \cite{survey}.

For trees, Kaplan, Lev, and Roditty~\cite{partition} proved that every tree with at most one vertex of degree two is 
antimagic \footnote{
Their original proof contains a minor error within a special case which was corrected by Liang, Wong, and Zhu~\cite{trees}.}  
In \cite{trees}, Liang, Wong, and Zhu~\cite{trees}  extended the results to trees with more vertices of degree 2. For instance it was proved that for a tree $T$ with $V_2(T)$ the set of degree-two vertices, if $V_2(T)$ induces a path, or both $V_2(T)$ and $V(T)-V_2(T)$ are independent sets, then $T$ is antimagic.

A {\em spider} is a tree with exactly one vertex of degree greater than two. It was proved  by Shang~\cite{Shang} and  Huang~\cite{Huang} independently that every spider is antimagic. A {\em double spider} is a tree with exactly two vertices of degree greater than two. Chang, Chin, Li, and  Pan~\cite{strongly} proved that every double spider is antimagic. A {\em caterpillar} is a tree so that removing all leaves resulting in a path. Deng and Y. Li \cite{caterpillar1} and Lozano, Mora, and Seara \cite{caterpillars} obtained results on caterpillars. 
Recently, Lozano, Mora, Seara, and Tey completely settled the problem that all caterpillars are antimagic~\cite{caterpillar all}.

Although \Href{connected antimagic} is widely open, there exist disconnected graphs that are not anitmagic.  It is easy to see that the disjoint union of two 3-paths $P_3$ is not antimagic. A {\em forest} is a graph where each component is a tree. When every component is a path, then it is called a {\em linear forest}. 
In~\cite{S18}, Shang proved that every linear forest without $P_2$ or $P_3$ as components is antimagic. A {\em star forest} is a graph where each component is a star.  Chen, Huang, Lin, Shang, and Lee~\cite{star forest} gave conditions of the components of a star forest to guarantee its antimagicness. 
Moreover, Dhananjaya~\cite{Dha} proved that the every double star forest is antimagic. 

\begin{definition}\rm
\label{shifted antimagic}
Let $G(V,E)$ be a graph with $m$ edges. Let $k$ be an integer, a {\em $k$-shifted antimagic labeling} is a bijection, $f: E  \to [k+1, k+m]$, such that the antimagic condition in \Href{antimagic} is satisfied. If $G$ admits a $k$-shifted antimagic labeling, then $G$ is {\em $k$-shifted antimagic}.  
\end{definition}

An antimagic labeling is a $k$-shifted antimagic labeling when $k=0$. The notion of {\em $k$-shifted antimagic labeling}  was introduced by Chang, Chen, Li, and Pan~\cite{shifted}, in which it was proved that every tree except $K_2$ is $k$-shifted antimagic for sufficiently large $k$.  

The aim of this article is to investigate $k$-shifted antimagic labelings for {\em spider forests} (each component is a spider).  
 For a spider, a path connecting its center and a leaf is called a {\em leg}. The length of a leg is the number of edges on that path. A length-one leg is called {\em 1-leg} (a pendant edge incident to the center). A star is a spider where each leg is a 1-leg.  As  stars are spiders, the results of Chen {\em et al.} on star forests~\cite{star forest} can be viewed as  special cases of our results when $k=0$. 

In this article, we prove the following three results: 

\begin{theorem}
\label{>2}
A spider forest where each spider has no  1-legs is $k$-shifted antimagic for all $k \geq 0$.
\end{theorem}

\begin{theorem}
\label{shiftedspider}
Let $G$ be a spider forest with $m$ edges. There exists a positive integer $k_0 < m$ so that $G$ is $k$-shifted antimagic for all $k\geq k_0$ and $k \leq -(m+k_0+1)$. \end{theorem}

\begin{theorem}
\label{1even}
Let $G$ be a spider forest such that the length of each leg of every spider is either 1 or even. Then $G$ is $k$-shifted antimagic for all $k \geq 0$.
\end{theorem}

The proofs of the above theorems are presented in the next section. More discussions and problems will be given in the last section.

\section{The Proofs}  


Throughout this section, we denote a spider forest by 
$ 
G=X_1 \cup X_2 \cdots \cup X_t
$  
where each $X_i$ is a spider centered at $w_i$. 
The number of edges of $G$ 
is denoted by $m$. 
An {\em odd leg} ({\em even leg}, respectively) is a leg of odd (even, respectively) length, and will be written as $v_1,v_2,\dots ,v_{2p_j+1},w_{i_j}$ ($u_1,u_2,\dots u_{2q_j}, w_{i_j}$, respectively).

%
\subsection{Proof of \Href{>2}}
%

Assume $G$ is a spider forest so that each $X_i$ contains legs of lengths at least two.  For each spider $X_i$ we reserve one of its legs, denoted by $L_i^*$, that will be labeled at the end of the process.  
If $X_i$ has an even leg, then we choose any  even leg as the reserved leg $L_i^*$ for $X_i$. Otherwise, $X_i$ contains only odd legs, we reserve any odd leg $L_i^*$ for $X_i$. 
Without loss of generality (re-arrange the spiders if needed) we may assume $L_i^*$ be odd legs for $i \in [1, t']$, and even legs for $i \in [t'+1, t]$. 

\begin{figure}[ht]
    \centering
\begin{picture}(360,100)
\textcolor{red}{
\put(120,80){\line(1,-1){20}}
\put(140,60){\line(0,-1){45}}
\put(30,80){\line(3,-2){30}}
\put(60,60){\line(0,-1){60}}
\put(180,60){\line(0,-1){15}}
\put(210,80){\line(-3,-2){30}}
\put(300,60){\line(0,-1){45}}
\put(320,80){\line(-1,-1){20}}
}
\put(30,90){$X_1$}
\put(120,90){$X_2$}
\put(210,90){$X_3$}
\put(320,90){$X_4$}
\put(30,80){\circle*{4}}
\put(30,80){\line(-1,-2){10}}
\put(30,80){\line(1,-2){10}}
\put(30,80){\line(-3,-2){30}}
\put(120,80){\circle*{4}}
\put(120,80){\line(0,-1){20}}
\put(120,80){\line(-1,-1){20}}
\put(210,80){\circle*{4}}
\put(210,80){\line(-1,-2){10}}
\put(210,80){\line(1,-2){10}}
\put(210,80){\line(3,-2){30}}
\put(320,80){\circle*{4}}
\put(320,80){\line(1,-1){20}}
\put(320,80){\line(2,-1){40}}
\put(320,80){\line(-2,-1){40}}
\put(320,80){\line(0,-1){20}}
\put(100,60){\circle*{4}}
\put(120,60){\circle*{4}}
\put(140,60){\circle*{4}}
\put(100,60){\line(0,-1){30}}
\put(120,60){\line(0,-1){30}}
\put(0,60){\circle*{4}}
\put(20,60){\circle*{4}}
\put(40,60){\circle*{4}}
\put(60,60){\circle*{4}}
\put(0,60){\line(0,-1){30}}
\put(20,60){\line(0,-1){30}}
\put(40,60){\line(0,-1){60}}
\put(180,60){\circle*{4}}
\put(200,60){\circle*{4}}
\put(220,60){\circle*{4}}
\put(240,60){\circle*{4}}
\put(200,60){\line(0,-1){30}}
\put(220,60){\line(0,-1){45}}
\put(240,60){\line(0,-1){60}}

\put(280,60){\circle*{4}}
\put(300,60){\circle*{4}}
\put(320,60){\circle*{4}}
\put(340,60){\circle*{4}}
\put(360,60){\circle*{4}}
\put(280,60){\line(0,-1){15}}
\put(320,60){\line(0,-1){45}}
\put(340,60){\line(0,-1){45}}
\put(360,60){\line(0,-1){60}}

\put(100,45){\circle*{4}}
\put(120,45){\circle*{4}}
\put(140,45){\circle*{4}}
\put(0,45){\circle*{4}}
\put(20,45){\circle*{4}}
\put(40,45){\circle*{4}}
\put(60,45){\circle*{4}}
\put(180,45){\circle*{4}}
\put(200,45){\circle*{4}}
\put(220,45){\circle*{4}}
\put(240,45){\circle*{4}}
\put(280,45){\circle*{4}}
\put(300,45){\circle*{4}}
\put(320,45){\circle*{4}}
\put(340,45){\circle*{4}}
\put(360,45){\circle*{4}}
\put(100,30){\circle*{4}}
\put(120,30){\circle*{4}}
\put(140,30){\circle*{4}}
\put(0,30){\circle*{4}}
\put(20,30){\circle*{4}}
\put(40,30){\circle*{4}}
\put(60,30){\circle*{4}}
\put(200,30){\circle*{4}}
\put(220,30){\circle*{4}}
\put(240,30){\circle*{4}}
\put(300,30){\circle*{4}}
\put(320,30){\circle*{4}}
\put(340,30){\circle*{4}}
\put(360,30){\circle*{4}}
\put(40,15){\circle*{4}}
\put(60,15){\circle*{4}}
\put(140,15){\circle*{4}}
\put(220,15){\circle*{4}}
\put(240,15){\circle*{4}}
\put(300,15){\circle*{4}}
\put(320,15){\circle*{4}}
\put(340,15){\circle*{4}}
\put(360,15){\circle*{4}}
\put(40,0){\circle*{4}}
\put(60,0){\circle*{4}}
\put(240,0){\circle*{4}}
\put(360,0){\circle*{4}}
\end{picture}
    \caption{A spider forest $G$ with a reserved leg $L_i^*$ (red) of each spider $X_i$.}
    \label{spiderforest}
\end{figure}
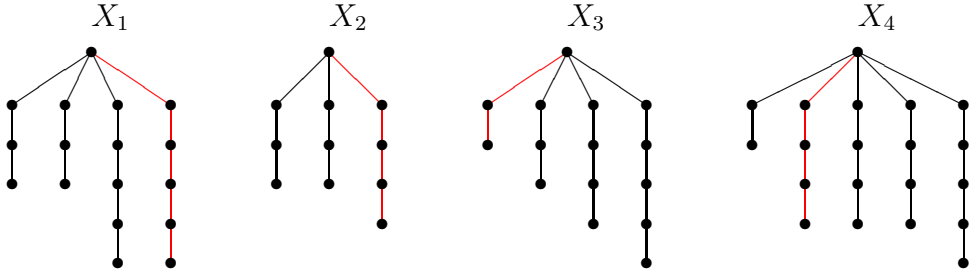

For the remaining non-reserved legs, we arrange them as $L_1, L_2,\ldots, L_d$, where $d=(\sum_{i=1}^t\deg(w_i))-t$, such that $L_i$ are odd legs for $i \in [1, n_1]$, and even legs for $i \in [n_1+1, d]$. Denote $\ell_i$ and $\ell_i^*$ the lengths of $L_i$ and $L_i^*$, respectively. See \Href{spiderforest} for an example. 

Denote 
$$
\begin{array}{llll}
a&=&\sum\limits_{i=1}^{n_1} \frac{\ell_i -1}{2} \\ 
b&=&\sum\limits_{i=n_1+1}^{d} \frac{\ell_i}{2}\\  
c_1&=&\sum\limits_{i=1}^{t'} \frac{\ell_i^*+1}{2}+\sum\limits_{i=t'+1}^{t} \frac{\ell_i^*}{2} \\
c_2&=&\sum\limits_{i=1}^{t'} \frac{\ell_i^* -1}{2} +\sum\limits_{i=t'+1}^{t} \frac{\ell_i^*}{2}.
\end{array}
$$

Let $k$ be a non-negative integer. 
In the following, we give a $k$-shifted antimagic labeling for $G$. Partition the set of labels to be used $[k+1, k+m]$ into sets of consecutive numbers, $[k+1, k+m]=I_1 \cup I_2\cup I_3\cup I_4\cup I_5\cup I_6\cup I_7$:
$$
\begin{array}{llll}
I_1&=&[k+1, k+a], &|I_1|=a\\   
I_2&=&[k+a+1, k+a+b], &|I_2|=b\\ 
I_3&=&[k+a+b+1,k+a+b+c_1],&|I_3|=c_1\\ 
I_4&=&[k+a+b+c_1+1, k+2a+b+c_1], &|I_4|=a \\  
I_5&=&[k+2a+b+c_1+1, k+2a+2b+c_1], &|I_5|=b\\
I_6&=&[k+2a+2b+c_1+1,k+2a+2b+c_1+c_2],&|I_6|=c_2\\  I_7&=&[k+2a+2b+c_1+c_2+1,k+m], &|I_7|=n_1.  
\end{array}
$$

\begin{figure}[ht]
     \centering
\begin{picture}(240,100)
\put(0,80){$\overbrace{\hspace{220pt}}^{\mbox{odd legs}}$}
\put(0,30){\circle*{4}}
\put(0,45){\circle*{4}}
\put(0,60){\circle*{4}}
\put(0,75){\circle*{4}}
\put(0,75){\line(0,-1){45}}
\put(30,30){\circle*{4}}
\put(30,45){\circle*{4}}
\put(30,60){\circle*{4}}
\put(30,75){\circle*{4}}
\put(30,75){\line(0,-1){45}}
\put(60,30){\circle*{4}}
\put(60,45){\circle*{4}}
\put(60,60){\circle*{4}}
\put(60,75){\circle*{4}}
\put(60,75){\line(0,-1){45}}
\put(90,30){\circle*{4}}
\put(90,45){\circle*{4}}
\put(90,60){\circle*{4}}
\put(90,75){\circle*{4}}
\put(90,75){\line(0,-1){45}}
\put(120,30){\circle*{4}}
\put(120,45){\circle*{4}}
\put(120,60){\circle*{4}}
\put(120,75){\circle*{4}}
\put(120,75){\line(0,-1){45}}
\put(150,0){\circle*{4}}
\put(150,15){\circle*{4}}
\put(150,30){\circle*{4}}
\put(150,45){\circle*{4}}
\put(150,60){\circle*{4}}
\put(150,75){\circle*{4}}
\put(150,75){\line(0,-1){75}}
\put(180,0){\circle*{4}}
\put(180,15){\circle*{4}}
\put(180,30){\circle*{4}}
\put(180,45){\circle*{4}}
\put(180,60){\circle*{4}}
\put(180,75){\circle*{4}}
\put(180,75){\line(0,-1){75}}
\put(210,0){\circle*{4}}
\put(210,15){\circle*{4}}
\put(210,30){\circle*{4}}
\put(210,45){\circle*{4}}
\put(210,60){\circle*{4}}
\put(210,75){\circle*{4}}
\put(210,75){\line(0,-1){75}}

\footnotesize{
\put(3,35){1}
\put(3,50){27}
\put(3,65){52}
\put(33,35){2}
\put(33,50){28}
\put(33,65){53}
\put(63,35){3}
\put(63,50){29}
\put(63,65){54}
\put(93,35){4}
\put(93,50){30}
\put(93,65){55}
\put(123,35){5}
\put(123,50){31}
\put(123,65){56}
\put(153,5){6}
\put(153,20){32}
\put(153,35){7}
\put(153,50){33}
\put(153,65){57}
\put(183,5){8}
\put(183,20){34}
\put(183,35){9}
\put(183,50){35}
\put(183,65){58}
\put(213,5){10}
\put(213,20){36}
\put(213,35){11}
\put(213,50){37}
\put(213,65){59}
}
\end{picture}
\quad
\begin{picture}(100,80)
\put(0,80){$\overbrace{\hspace{100pt}}^{\mbox{even legs}}$}
\put(0,45){\circle*{4}}
\put(0,60){\circle*{4}}
\put(0,75){\circle*{4}}
\put(0,75){\line(0,-1){30}}
\put(30,15){\circle*{4}}
\put(30,30){\circle*{4}}
\put(30,45){\circle*{4}}
\put(30,60){\circle*{4}}
\put(30,75){\circle*{4}}
\put(30,75){\line(0,-1){60}}
\put(60,15){\circle*{4}}
\put(60,30){\circle*{4}}
\put(60,45){\circle*{4}}
\put(60,60){\circle*{4}}
\put(60,75){\circle*{4}}
\put(60,75){\line(0,-1){60}}
\put(90,15){\circle*{4}}
\put(90,30){\circle*{4}}
\put(90,45){\circle*{4}}
\put(90,60){\circle*{4}}
\put(90,75){\circle*{4}}
\put(90,75){\line(0,-1){60}}
\footnotesize{
\put(3,50){12}
\put(3,65){38}
\put(33,20){13}
\put(33,35){39}
\put(33,50){14}
\put(33,65){40}
\put(63,20){15}
\put(63,35){41}
\put(63,50){16}
\put(63,65){42}
\put(93,20){17}
\put(93,35){43}
\put(93,50){18}
\put(93,65){44}
}
\end{picture}
     \caption{Labeling the non-reserved odd and even legs of the spider forest in \Href{spiderforest} with labels in $[1, 59] \setminus ([19,26] \cup [45,51])= I_1 \cup I_2 \cup I_4 \cup I_5 \cup I_7 =[1,11] \cup [12,18] \cup [27,37] \cup [38,44] \cup [52,59]$.} 
     \label{oddevenlegs}
 \end{figure}
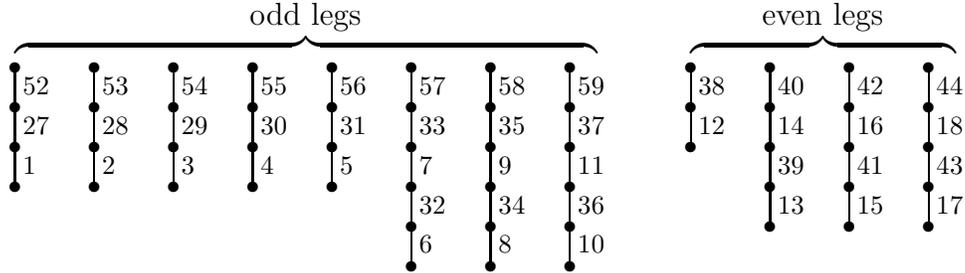

\noindent
Label the odd legs $L_1, L_2,\ldots,L_{n_1}$ using integers in $I_1 \cup I_4 \cup I_7$ by the following three steps: 
\begin{enumerate}[{\bf Step 1.}]
    \item Label $v_1v_2, v_3v_4, \ldots$, $v_{2p_j-1}v_{2p_j}$ of each path one by one consecutively (one path after another) with integers in $I_1$ consecutively in an increasing order. \item Label edges $v_2v_3, v_4v_5, \ldots,  v_{2p_j}v_{2p_j+1}$ with the integers in $I_4$ following the same order as in Step 1. 
    \item Label $v_{2p_j+1}w_{i_j}$ with $k+m-n_1+j$ in $I_7$ for  $ j \in[1, n_1]$.
\end{enumerate}
Similarly, we label the even legs $L_j$, $j \in [n_1+1, d]$, with numbers in $I_2$ and $I_5$ by Steps 1 and 2 in the above, regarding $v_{2p_j+1}$ as $w_{j}$.  
See \Href{oddevenlegs} for an example. 
 
By now, all edges are labelled, except the ones on the reserved leg $L_i^*$ for each $X_i$. Denote $\phi'(w_i)$ the partial sum of the labels of edges incident to $w_i$, $i \in [1,t]$. We label the edges of $L_i^*$, $i \in [1, t]$, 
by an iterated process using labels in $I_3 \cup I_6$. In each iteration, we complete the labeling for one spider. 

Initially, for each spider $X_i$, $i \in [1,t]$, define  
\[
\phi_1(w_i)=
\left\{
\begin{array}{ll}
\phi'(w_i)+k+a+b+c_1, &\mbox{ if }L_i^* \mbox{ is odd},  \\
\phi'(w_i)+k+2a+2b+c_1+c_2, & 
\mbox{ if }L_i^* \mbox{ is even}.
\end{array}
\right.
\]
Note that $k+a+b+c_1$ and 
$k+2a+2b+c_1+c_2$ are the maximum numbers in $I_3$ and $I_6$,   respectively. Denote $Z_1=\{X_i \ | \ \phi_1(w_i) \geq \phi_1(w_j) \  {\rm for \ all} \   j \in [1, t]\}$.  

{\bf Case 1.} Assume $L^*_i$ is odd for every $X_i \in Z_1$. Then we choose one $X_{i_1} \in Z_1$, and label the edge $e$ of $L^*_{i_1}$ incident to $w_{i_1}$ by $k+a+b+c_1$. Next we use Steps 1 and 2 in the above to label the remaining edges in $L^*_{i_1}$ by using the largest unassigned labels in $I_3$ and $I_6$, respectively. Note that the label $k+a+b+c_1+c_2$ is assigned to the edge incident to $e$. Hence, $X_{i_1}$ is completely  labeled. 

{\bf Case 2.} Assume there exists some $X_i \in Z_1$ where $L_{i}^*$ is even. Let $i_1=i$. Then label $L_{i_1}^*$ by Steps 1 and 2, using the largest unused labels in $I_3$ and $I_6$, respectively.  Note that, the label $k+2a+2b+c_1+c_2$ is assigned to the edge $e$ of $L_{i_1}^*$ incident to $w_{i_1}$, and $k+a+b+c_1$ is assigned to the edge incident to $e$ in $L_{i_1}^*$. 

\begin{figure}[ht]
    \centering
\begin{picture}(360,90)
\textcolor{red}{
\put(120,80){\line(1,-1){20}}
\put(140,60){\line(0,-1){45}}
\put(30,80){\line(3,-2){30}}
\put(60,60){\line(0,-1){60}}
\put(180,60){\line(0,-1){15}}
\put(210,80){\line(-3,-2){30}}
\put(300,60){\line(0,-1){45}}
\put(320,80){\line(-1,-1){20}}
}
\put(30,90){$X_1$}
\put(120,90){$X_2$}
\put(210,90){$X_3$}
\put(320,90){$X_4$}
\put(30,80){\circle*{4}}
\put(30,80){\line(-1,-2){10}}
\put(30,80){\line(1,-2){10}}
\put(30,80){\line(-3,-2){30}}
\put(120,80){\circle*{4}}
\put(120,80){\line(0,-1){20}}
\put(120,80){\line(-1,-1){20}}
\put(210,80){\circle*{4}}
\put(210,80){\line(-1,-2){10}}
\put(210,80){\line(1,-2){10}}
\put(210,80){\line(3,-2){30}}
\put(320,80){\circle*{4}}
\put(320,80){\line(1,-1){20}}
\put(320,80){\line(2,-1){40}}
\put(320,80){\line(-2,-1){40}}
\put(320,80){\line(0,-1){20}}
\put(100,60){\circle*{4}}
\put(120,60){\circle*{4}}
\put(140,60){\circle*{4}}
\put(100,60){\line(0,-1){30}}
\put(120,60){\line(0,-1){30}}
\put(0,60){\circle*{4}}
\put(20,60){\circle*{4}}
\put(40,60){\circle*{4}}
\put(60,60){\circle*{4}}
\put(0,60){\line(0,-1){30}}
\put(20,60){\line(0,-1){30}}
\put(40,60){\line(0,-1){60}}
\put(180,60){\circle*{4}}
\put(200,60){\circle*{4}}
\put(220,60){\circle*{4}}
\put(240,60){\circle*{4}}
\put(200,60){\line(0,-1){30}}
\put(220,60){\line(0,-1){45}}
\put(240,60){\line(0,-1){60}}

\put(280,60){\circle*{4}}
\put(300,60){\circle*{4}}
\put(320,60){\circle*{4}}
\put(340,60){\circle*{4}}
\put(360,60){\circle*{4}}
\put(280,60){\line(0,-1){15}}
\put(320,60){\line(0,-1){45}}
\put(340,60){\line(0,-1){45}}
\put(360,60){\line(0,-1){60}}

\put(100,45){\circle*{4}}
\put(120,45){\circle*{4}}
\put(140,45){\circle*{4}}
\put(0,45){\circle*{4}}
\put(20,45){\circle*{4}}
\put(40,45){\circle*{4}}
\put(60,45){\circle*{4}}
\put(180,45){\circle*{4}}
\put(200,45){\circle*{4}}
\put(220,45){\circle*{4}}
\put(240,45){\circle*{4}}
\put(280,45){\circle*{4}}
\put(300,45){\circle*{4}}
\put(320,45){\circle*{4}}
\put(340,45){\circle*{4}}
\put(360,45){\circle*{4}}
\put(100,30){\circle*{4}}
\put(120,30){\circle*{4}}
\put(140,30){\circle*{4}}
\put(0,30){\circle*{4}}
\put(20,30){\circle*{4}}
\put(40,30){\circle*{4}}
\put(60,30){\circle*{4}}
\put(200,30){\circle*{4}}
\put(220,30){\circle*{4}}
\put(240,30){\circle*{4}}
\put(300,30){\circle*{4}}
\put(320,30){\circle*{4}}
\put(340,30){\circle*{4}}
\put(360,30){\circle*{4}}
\put(40,15){\circle*{4}}
\put(60,15){\circle*{4}}
\put(140,15){\circle*{4}}
\put(220,15){\circle*{4}}
\put(240,15){\circle*{4}}
\put(300,15){\circle*{4}}
\put(320,15){\circle*{4}}
\put(340,15){\circle*{4}}
\put(360,15){\circle*{4}}
\put(40,0){\circle*{4}}
\put(60,0){\circle*{4}}
\put(240,0){\circle*{4}}
\put(360,0){\circle*{4}}
\tiny
{
\put(3,35){1}
\put(3,50){27}
\put(2,65){52}
\put(23,35){2}
\put(23,50){28}
\put(23,65){53}
\put(143,20){19}
\put(143,35){45}
\put(143,50){20}
\put(137,65){46}
\put(103,35){3}
\put(103,50){29}
\put(107,65){54}
\put(123,35){4}
\put(123,50){30}
\put(123,65){55}
\put(43,5){6}
\put(43,20){32}
\put(43,35){7}
\put(43,50){33}
\put(40,65){57}
\put(63,5){21}
\put(63,20){47}
\put(63,35){22}
\put(63,50){48}
\put(55,65){23}
\put(183,50){24}
\put(183,65){49}
\put(203,35){5}
\put(203,50){31}
\put(203,65){56}
\put(243,5){8}
\put(243,20){34}
\put(243,35){9}
\put(243,50){35}
\put(235,65){58}
\put(363,5){10}
\put(363,20){36}
\put(363,35){11}
\put(363,50){37}
\put(355,65){59}
\put(343,20){17}
\put(343,35){43}
\put(343,50){18}
\put(335,65){44}
\put(323,20){15}
\put(323,35){41}
\put(323,50){16}
\put(323,65){42}
\put(303,20){25}
\put(303,35){50}
\put(303,50){26}
\put(308,65){51}
\put(283,50){12}
\put(283,65){38}
\put(223,20){13}
\put(223,35){39}
\put(223,50){14}
\put(220,65){40}
}
\end{picture}
    \caption{At the end we obtain an antimagic labeling for the spider forest in \cref{spiderforest}}
    \label{labelingThm1}
\end{figure}
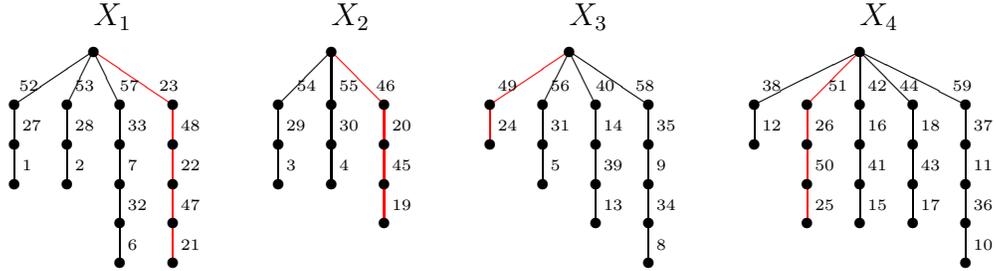

In both cases, $\phi(w_{i_1})=\phi_1(w_{i_1})$. 
In the next iteration, similar to the initial iteration, for the remaining spiders we define $\phi_2(w_i)$ to be $\phi'(w_i)$ plus the maximum unused labels in $I_3$ and $I_6$,  according to the parity of $L_i^*$, for all $j \neq i_1$.  Then we choose a spider $X_{i_2}$ with the maximum value of $\phi_2(w_{i_2})$ and fix the label of the edge incident to  $w_{i_2}$ (then all the edges on $L_{i_2}^*$)  and so  
$\phi(w_{i_1})$ is determined. 
Note that $\phi(w_{i_2}) < \phi(w_{i_1})$. After $t$ iterations, $\phi(w_i)$ are  determined for all $i \in [1,t]$ and they are all distinct. 
See \Href{labelingThm1}.

It remains to show that the labeling satisfies the antimagic condition. 
Each pendent vertex receives a label from $I_1\cup I_2\cup I_3$ uniquely as its vertex-sum, and a non-leaf vertex is adjacent to some edge whose label is in $I_4\cup I_5\cup I_6\cup I_7$.
So the vertex-sum of a pendent vertex is smaller than that of any non-leaf vertex. 
Next, we claim that degree-2 vertices all have distinct vertex-sums. For every degree-2 vertex $v$, denote the two edges incident to $v$ by $e_v$ and $e_{v'}$ where $f(e_v) < f(e'_{v})$. For two degree-2 vertices $u \neq v$, if $\{f(e_v), f(e'_{v})\} \cap \{f(e_u), f(e'_{u})\} \neq \emptyset$ then clearly $\phi(v) \neq \phi(u)$. Assume $\{f(e_v), f(e'_{v})\} \cap  \{f(e_u), f(e'_{u})\} = \emptyset$.  
It is easy to see from the ordering of our labeling, it holds that $f(e_v) < f(e_u)$ if and only if $f(e'_{v}) < f(e'_{u})$. Hence all degree-2 vertices have distinct vertex-sums.

Now we claim that $\phi(v) < \phi(w_i)$ holds for any degree-2 vertex $v$ and $i \in [1,t]$. 
Note that $\phi(v)$ is at most  
the sum of the maximum in $I_4$ and the maximum in $I_7$. That is, $\phi(v) \leq (k+2a+b+c_1)+(k+m)=m+2k+2a+b+c_1$.  
The lower bound of $\phi(w_i)$ of $X_i$ can be estimated by assuming that $X_i$ has only 3 legs. 
Let $j_1 < j_2 < j_3$ be the three labels of edges incident to $w_i$. 
If $X_i$ has three odd legs, then $j_1 \in I_3$ and $j_2, j_3 \in I_7$. 
If $X_i$ has one odd and two even legs, then $j_1 \in I_5$, $j_2 \in I_6$, $j_3 \in I_7$. 
If $X_i$ has two odd and one even legs, then $j_1 \in I_6$ and $j_2, j_3 \in I_7$. If $X_i$ has three even legs, then $j_1, j_2 \in I_5$ and $j_3 \in I_6$. Direct calculation shows that in each possibility $\phi(w_i) > m+2k+2a+b + c_1$. Therefore, the proof of \Href{>2} is complete.

%
\subsection{Proof of \Href{shiftedspider}}
%

The proof of \Href{shiftedspider} is similar to the one of \Href{>2}. 
Assume $G$ is a spider forest with $m$ edges. 
For each $i \in [1,t]$, suppose $X_i$ has $s_i$ 1-legs and $r_i$ legs of length at least two.  
Denote  $s=\sum_{i=1}^ts_i$. 
Let $\alpha$, $\beta$ and $\gamma$ be the numbers of spiders $X_i$ with $r_i=0$, $r_i=1$ and $r_i=2$, respectively.

For each spider $X_i$, we reserve one of its legs $L^*_i$ that will be labeled at the end of the labeling process.  For an $X_i$ with $r_i\leq 2$, we choose a 1-leg as the reserved leg $L^*_i$. For $X_i$ with $r_i \geq 3$, if $X_i$ has an even leg, then we choose an even leg as the reserved leg; otherwise, we choose an odd leg as the reserved leg.  By re-arranging the spiders, we may assume $L^*_1, L^*_2, \ldots , L^*_{t'}$ are odd legs  (including the $\alpha+\beta+\gamma$ 1-legs) and $L^*_i$ are  even legs for $i \in [t'+1, t]$. 

Additionally, for a spider $X_i$ with $r_i=0$, we collect two non-reserved 1-legs; and for $r_i=1$, we collect one non-reserved 1-leg. Arrange these collected legs in the order  $L_1, L_2, \ldots, L_{2\alpha+\beta}$ such that $L_i$ and $ L_{i+\alpha}$ are from the same spider with $r_i=0$,  $i\in[1, \alpha]$.

Arrange the remaining non-reserved and non-collected legs in the order   
$L_{2\alpha+\beta+1}, \ldots, L_d$, where $d=(\sum^t_{i=1}\deg(w_i))-t$, such that  $L_i$ are odd legs for $i\in[1, n_1]$ and particularly  $L_1,L_2,\ldots, L_{s-(\alpha+\beta+\gamma)}$ are the 1-legs. Denote $\ell_i$ and $\ell^*_i$ the lengths of legs $L_i$ and $L^*_i$, 
respectively.  Denote 
$$
\begin{array}{llll} 
a&=&\sum\limits_{i=1}^{n_1} \frac{\ell_i -1}{2} \\
b&=&\sum\limits_{i=n_1+1}^{d} \frac{\ell_i}{2} \\
q&=&s-(3\alpha+2\beta+\gamma)\\ c_1&=&\sum\limits_{i=1}^{t'} \frac{\ell_i^*+1}{2}  +\sum\limits_{i=t'+1}^{t} \frac{\ell_i^*}{2}  \\
c_2&=&\sum\limits_{i=1}^{t'} \frac{\ell_i^*-1}{2}  +\sum\limits_{i=t'+1}^{t} \frac{\ell_i^*}{2}. 
\end{array} 
$$

Let $k$ be a positive integer, $k \geq n_1+a+c_2-2q$. It is easy to see that $n_1+a+c_2-2q < m$. 
In the following, we give a $k$-shifted antimagic labeling for $G$. Partition the set of labels to be used $[k+1,k+m]$ into sets of consecutive numbers, $I_1\cup I_2 \cup \ldots \cup I_9$ by:

$$
\begin{array}{lll}
I_1&=&[k+1, k+a], \\   
I_2&=&[k+a+1, k+a+b], \\ 
I_3&=&[k+a+b+1,k+a+b+q],\\ 
I_4&=&[k+a+b+q+1, k+a+b+q+c_1],\\  
I_5&=&[k+a+b+q+c_1+1, k+a+b+q+c_1+(2\alpha+\beta)],\\
I_6&=&[k+a+b+q+c_1+(2\alpha+\beta)+1,k+2a+b+q+c_1+(2\alpha+\beta)],\\ 
I_7&=&[k+2a+b+q+c_1+(2\alpha+\beta)+1, k+2a+2b+q+c_1+(2\alpha+\beta)],\\
I_8&=&[k+2a+2b+q+c_1+(2\alpha+\beta)+1,k+2a+2b+q+c_1+(2\alpha+\beta)+c_2],\\  I_9&=&[k+2a+2b+q+c_1+(2\alpha+\beta)+c_2+1,k+m]. 
\end{array}
$$
Observe that $|I_1|=|I_6|=a$, $|I_2|=|I_7|=b$, $|I_3|=q$, $|I_4|=c_1$, $|I_5|=2\alpha+\beta$, $|I_8|=c_2$, and $|I_9|=n_1-q-(2\alpha+\beta)$. 

We first label the non-reserved 1-legs using numbers from $I_3 \cup I_5$ as follows. For $i \in [1, \alpha]$, label $L_i$ and $L_{\alpha+i}$ with $k+a+b+q+c_1+i$ and $k+a+b+q+c_1+(2\alpha+\beta)+1-i$, respectively; for $i \in [1, \beta]$,  label $L_{2\alpha+i}$ with  $k+a+b+q+c_1+\alpha+i$; and  
label each of the remaining 1-leg   
$L_i$, $i \in [2\alpha+\beta+1, s-(\alpha+\beta+\gamma)]$, with a number in $I_3$.

For the remaining odd legs  
$L_i$ of length at least three, $i \in [s-(\alpha+\beta+\gamma)+1, n_1]$, we label them the same as Steps 1, 2, 3  in the proof of \Href{>2}, using the numbers from $I_1$, $I_6$ and $I_9$, respectively.  Similarly, we label the remaining even legs $L_i$, $i \in [n_1+1, d]$, by Steps 1 and 2 in \Href{>2} (treating $v_{2p_i+1}$ as $w_i$),   using labels in $I_2$ and $I_7$, respectively.

Thus far we have labeled all edges, except the reserved legs 
$L_i^*$ of each spider $X_i$, $i \in [1,t]$. We apply the same iterated process as in the proof of \Href{>2} to complete the labeling, using labels in $I_4$ and $I_8$, so that  $\phi(w_i)$ are all distinct for $i \in [1,t]$. Note that by the same comparison argument in the proof of~\Href{>2}, we see the vertex-sums of all pendent vertices and degree-2 vertices are distinct as well.

It remains to show that $\phi(w_i) >  \phi(v)$ for any degree-2 vertex $v$ and for any $i \in [1,t]$. The maximum $\phi(v)$ occurs when $v$ is incident to edges with a label in $I_6$ and the other in $I_9$. That is, 
\[
\phi(v) \leq  (k+2a+b+q+c_1+(2\alpha+\beta))+(k+m).
\]
The lower bound of $\phi(w_i)$ can be estimated by assuming some $X_i$ is $S_3$.
Thus  $\phi(w_i) \geq (k+a+b+q+1)+(k+a+b+q+c_1+j)+(k+a+b+q+c_1+(2\alpha+\beta)+1-j)$. By our assumption that $k\geq n_1+a+c_2-2q$ and since \begin{eqnarray*}
m&=&2a+2b+q+c_1+(2\alpha+\beta)+c_2+(n_1-q-(2\alpha+\beta))\\
&=&2a+2b+c_1+c_2+n_1,
\end{eqnarray*}
we have $\phi(w_i)>  (k+2a+b+q+c_1+(2\alpha+\beta))+(k+m)\geq \phi(v)$. 


The above $k$-shifted antimagic labeling $f$ for $G$ induces a $(-k-m-1)$-shifted antimagic labeling $f': E(G) \to [-(k+m), -(k+1)]$ by $f'(e) = -f(e)$ for every $e \in E(G)$.   Thus, we conclude that $G$ is $k$-shifted antimagic for $k \geq k_0$ or $k\leq -(m+k_0+1)$. This completes the proof of \Href{shiftedspider}.

%
\subsection{Proof of \Href{1even}}
%

Assume $G$ is a spider forest with $m$ edges, where the length of each leg is either 1 or even.  
Suppose each spider $X_i$ has $s_i$ 1-legs. Denote the total number of 1-legs in $G$ by $s=\sum_{i=1}^t s_i$. By the assumption of $G$, 
$m$ and $s$ have the same parity. 

Classify the spiders into four types:  
An $S_3$ spider is {\em Type-A}, an $S_n$ with $n \geq 4$ is {\em Type-B}, 
a non-star spider $X_i$ with $s_i \geq 2$ is {\em Type-C}, and  all other spiders are {\em Type-D}. Note that each Type-D spider has $s_i \leq 1$, so it contains at least two even legs. Without loss of generality, we may assume that  
$X_i$ is Type-A, Type-B, Type-C, and Type-D for 
$i$ in  $[1, t_1]$, $[t_1 + 1, t_2]$, $ [t_2+1, t_3]$ and $[t_3+1,t]$, respectively.

 Let $k$ be a non-negative integer. 
In the following, we give a $k$-shifted antimagic labeling for $G$. Partition the set of labels to be used $[k+1, k+m]$ into three sets of consecutive numbers, $[k+1,k+m]=I_1 \cup I_2\cup I_3$, where 
$$
\begin{array}{llll}
I_1 &=&[k+1, k+ \frac{m-s}{2}], &|I_1|= \frac{m-s}{2}\\  I_2 &=&[k+ \frac{m-s}{2}+1, k+\frac{m+s}{2}], &|I_2|=s\\    I_3 &=&[k+\frac{m+s}{2}+1, k+m],  &|I_3|=\frac{m-s}{2}.
\end{array} 
$$

For each Type-A or -B spider (a star) $X_i$, $i \in [1,t_2]$, we label two edges with numbers from $I_2$, $k+\frac{m-s}{2}+i$ and $k+\frac{m+s}{2}+1-i$. For a Type-A spider  $X_i$, $i \in [1,t_1]$, we label the third 1-leg by $k+\frac{m-s}{2}+t_2+i$.  By now, all Type-A spiders are completely labeled. The remaining unused labels in $I_2$ are $[k+\frac{m-s}{2}+t_2+t_1+1, k+\frac{m+s}{2}-t_2]$. 

For each Type-B or -C spider, we reserve two unlabeled  
1-legs that will be labeled at the end of the process. 
For each Type-D spider $X_j$, we reserve two even legs and denote them as  $L_1^{(j)}$ and $L_2^{(j)}$.  

Now arrange all remaining non-reserved even legs in an arbitrary order $L_1, L_2, \ldots$. Label these legs by the same method in Steps 1 and 2 in \Href{>2} (treating $v_{2p_i+1}$ as $w_j$) starting from  the smallest labels in $I_1$ and $I_3$, respectively. Afterwards, the remaining unused labels in $I_1$  (in $I_3$, respectively) form a set of consecutive integers $I_1'$ (or $I_3'$, respectively),
and $|I_1'| = |I_3'|$.

Next, label the non-reserved 1-legs with the smallest  available integers in $I_2$, starting from 
$k+\frac{m-s}{2}+t_1+t_2+1$ (the the smallest unused numbers in $I_2$).  

Thus far, all edges incident to  $w_i$ for $i \in [t_1+1, t]$ are labeled, except the two reserved legs for each $X_i$.  Denote the sum of these known labels incident to $w_i$ as $\phi'(w_i)$. By re-ordering the indices of spiders if necessary, we may assume that $X_i$ for $i \in [t_1+1, t_3]$ are type-B or -C, while $i \in [t_3+1, t]$ are type-D, where \[
\phi'(w_{t_1+1})\leqslant\phi'(w_{t_1+2})\leqslant
\cdots\leqslant\phi'(w_{t_3}) \ \ \mbox{ and } \ \ 
\phi'(w_{t_3+1})\leqslant
\phi'(w_{t_3+2})\leqslant
\cdots\leqslant\phi'(w_t). \eqno(*)
\]
Now we label the reserved even legs $L_1^{(t_3+1)},L_2^{(t_3+1)},L_1^{(t_3+2)},L_2^{(t_3+2)},\ldots,L_2^{(t)}$ of Type-D spiders $X_{t_3+1},\ldots, X_t$,  using Steps 1 and 2 in the proof of \Href{>2} (regarding $v_{2p_j+1}$ as $w_j$) 
with labels in $I_1'$ and $I_3'$,  respectively.

Denote $m'= k+\frac{m+s}{2}-t_2-2t_3$. 
For a Type-B or -C spider $X_i$, $i \in[t_1+1, t_3]$, we temporarily label the two unlabeled 1-legs with $m'+2i-1$ and $m'+2i$, from $I_2$.  
Now all edges incident to the centers are labeled, and the vertex-sums satisfy the following:

\begin{description}
\item{(1)} 
$\phi(w_{i+1})-\phi(w_i)=1$ for   $i \in [1, t_1-1]$,
\item{(2)} 
$\phi(w_{i+1})-\phi(w_i)\geq 3$ for $i=t_1$,
\item{(3)} 
$\phi(w_{i+1})-\phi(w_i)\geq 4$ for 
$i \in [t_1+1, t_3-1] \cup [t_3+1, t-1]$.
\end{description}

From our labeling scheme, (1) and (2) are obvious. 
Inequality in (3)  holds because for each $i$ the last two labels assigned to edges incident to $w_i$
are less than those to $w_{i+1}$.  

In the following, we prove $\phi(w_i) \neq \phi(w_j)$ if $i \neq j$. If $X_i$ is Type-A, then $\phi(w_i)$ is the sum of three numbers from $I_2$. The center of a Type-D spider receives at least two labels from $I_3$,  and even if it receives one label from $I_2$, the label remains to be greater than at least two labels of each Type-A spider.
Hence $\phi(w_i)<\phi(w_j)$, if $X_i$ is Type-A and $X_j$ is Type-D. 

Assume  $\phi(w_i)=\phi(w_j)$ for some $i\neq j$.  By (*),  
one of them must be Type-B or -C and the other is Type-D. We call a Type-B or -C spider $X_i$ {\em trouble} if $\phi(w_i)=\phi(w_j)$ for some Type-D spider $X_j$. A non-trouble spider is called {\it good}. 

We make adjustments so that all Type-B or -C spiders become good.  Let $X_i$ be a trouble spider with the  smallest index $i$. 
We switch the special labels $m'+2i$ in $X_i$ with the label $m'+2i+1$ in $X_{i+1}$. That is, the vertex-sum of $w_{i}$ is increased by 1, while the vertex-sum of $w_{i+1}$ is  decreased by 1. By (3), $X_i$ becomes good after the switch. 
Continue this process until  $X_{t_3-2}$. 

It suffices to consider the case that all spiders $X_q$ are good for $t_1+1 \leq q \leq t_3-2$. 
If both $X_{t_3-1}$ and $X_{t_3}$ are good, we are done.  Assume at least one of them is a trouble spider.   
Denote the current updated vertex-sums for $w_{t_3-1}$ and $w_{t_3}$ by $\phi(w_{t_3-1})$ and $\phi(w_{t_3})$, respectively.   
Note that $m'+2t_3-2$ is a special label for $X_{t_3-1}$,  while $m'+2t_3$ and $m'+2t_3-1$ are the two  special labels for $w_{t_3}$. Similar to the above, we switch the special labels of $m'+2t_3-2$ and $m'+2t_3-1$ between $X_{t_3-1}$ and $X_{t_3}$. Then the new vertex-sums for $X_{t_3-1}$ and $X_{t_3}$ are $\phi(w_{t_3-1})+1$ and $\phi(w_{t_3})-1$, respectively.  
If $X_{t_3-1}$ or $X_{t_3}$ was a trouble spider, then it  becomes good after the switch.  

Hence it remains to consider the situation that $X_{t_3-1}$ was trouble before the switch, while $X_{t_3}$ becomes trouble after the switch (the other case that $X_{t_3}$ was trouble, while $X_{t_3-1}$ becomes trouble after the switch can be proved similarly). That is,  $\phi(w_{t_3-1})=\phi(w_{u'})$ and $\phi(w_{t_3}) -1 = \phi(w_{u})$ for some $t_3 +1\leq u' < u \leq t$. By (3),  $\phi(w_{u'}) \leq \phi(w_{u})-4$, so we get   $\phi(w_{t_3-1}) \leq \phi(w_{t_3})-5$. For this case, instead of switching the special labels $m'+2t_3-1$ with $m'+2t_3-2$, we switch $m'+2t_3$ with $m'+2t_3-2$ between $X_{t_3}$ and $X_{t_3-1}$.  This implies that the new vertex-sums of $w_{t_3-1}$ and $w_{t_3}$ are $\phi(w_{t_3-1})+2$ and $\phi(w_{t_3})-2$, respectively. Afterwords both $X_{t_3-1}$ and $X_{t_3}$ become good. Further, as $\phi(w_{t_3-1}) \leq \phi(w_{t_3})-5$, it is clear that $w_{t_3-1}$ and $w_{t_3}$ have distinct new vertex-sums. Hence, after the process, we conclude that $\phi(w_i) \neq \phi(w_j)$ for $i \neq j$.

We claim that the above process gives a $k$-shifted antimagic labeling. It is obvious that the vertex-sums of the leaves are all distinct. 
For a degree-2 vertex $u_i$ in an even leg $u_1u_2\ldots u_{2q_j}w_{i_j}$, if $i$ is even, then $\phi(u_i)=2\ell+ \frac{m+s}{2}$ for some $\ell$; if $i$ is odd, then  $\phi(u_i)=2\ell'+ \frac{m+s}{2}-1$ for some $\ell'$.  The values of $\ell$ and $\ell'$ are uniquely determined by $u_i$ and in $I_1$. 
By direct calculation, one can verify that they are all distinct and greater than $k+\frac{m+s}{2}$.   
As the vertex-sum of a leaf is at most $k+\frac{m+s}{2}$, the the vertex-sum of a degree-2 vertex is greater than that of a leaf. 

It remains to show that  $\phi(w_i) > \phi(v)$ for any degree-2 vertex $v$ and for any $i \in [1,t]$. Note that $\phi(v) \leq (k+\frac{m-s}{2})+(k+m)$, the sum of the largest integers in $I_1$ and $I_3$.  
The value of $\phi(w_i)$ is at least  $(k+\frac{m-s}{2}+i)+(k+\frac{m+s}{2}+1-i)+(k+\frac{m-s}{2}+t_1+i)$ if $X_i$ is Type-A or -B, 
or it is the sum of at least two labels in $I_3$ and one label in $I_2 \cup I_3$ if $X_i$ is Type-C or -D. In all cases, 
it is easy to see that $\phi(w_i) > \phi(v)$, since the sum of two labels $I_3$ is greater than $m+k$. Therefore, the proof for \Href{1even} is complete.

\section{Discussion and Open Problems}
%

A basic idea behind our labeling schemes is to classify the labels as small, middle, and large. Each pendant edge is labeled with a small label, while the edges incident to a degree-2 vertex receive a small and a middle labels, or a middle and a large labels, and the edges incident to the  center of a spider are labeled with  middle or large labels.

When applying the above idea to label the stars in a spider forest   difficulties arose as every leg of a star is both a pendant edge and an edge incident to the center. Our original intent might lead to a small vertex-sum of the center of a star, increasing the possibility for it to be the same as the vertex-sum of a degree-2 vertex. One way to resolve this problem is to   label the stars with the middle or large labels, giving larger vertex-sums of the centers of the stars. However, this might create the problem that the vertex-sums of the leaves of the stars might coincide with the vertex-sums of some degree-2 vertices of other non-star spiders.  

 To resolve both issues, we use larger  labels so that the vertex-sum of any center (with degree at least 3) will always be large enough.   For example, if we use labels in $[2m+1,3m]$, then the vertex-sum of a center is always greater than that of a degree-2 or degree-1 vertex.  This motivates us to shift the labels from $[1,m]$ to $[k+1,k+m]$,  and then minimize the values of $k$.  \Href{shiftedspider} for general spider forests came out from this approach. 

In addition to the conditions in  \Href{>2} and \Href{1even}, we could 
impose other conditions on the numbers of odd legs and even legs of the spiders, and show those spider forests 
are antimagic. However, the conditions become relatively complicated and tedious, and the arguments are analogously lengthy. 
Nevertheless, 
We believe that every spider forest is antimagic and propose the following conjecture.

\begin{conjecture}\label{conj}
Every spider forest is $k$-shifted antimagic for all $k\geq 0$. 
\end{conjecture}

Recall that a star forest is also a  spider forest. Recently, Dhananjaya and Li~\cite{consecutive} proved that if a star forest is not one of the following cases: a star, or the union of two $S_3$'s, or the union of three $S_3$'s, or the union of two $S_4$'s, and every component is not an $S_2$, then it is $k$-shifted antimagic for every integer $k$. 
For spiders other than stars, an example in~\cite{shifted} shows that a spider with one  length-1 leg and two length-2 legs is $k$-shifted antimagic if and only if $k\neq -3$.  
A more challenging problem than  Conjecture~\ref{conj} 
might be to  determine all the integers $k$ for which a given spider forest is $k$-shifted antimagic.

\end{document}